\documentclass[12pt]{amsart}

\usepackage{amssymb}
\def\Z{{\mathbb Z}}
\def\Aut{{\rm Aut}}
\def\PSL{{\rm PSL}}

\newcommand{\Out}{\mathop{\textrm{Out}}}

\def\PGamL{{\rm P\Gamma L}}
\newcommand{\Sym}{\mathop{\textrm{Sym}}}
\newcommand{\Inn}{\mathop{\textrm{Inn}}}
\def\normal{\triangleleft}
\theoremstyle{definition}

\numberwithin{equation}{section}

\newtheorem{thm}[equation]{Theorem}
\newtheorem{defn}[equation]{Definition}
\newtheorem{lem}[equation]{Lemma}
\newtheorem{cor}[equation]{Corollary}
\newtheorem{prop}[equation]{Proposition}
\newtheorem{hey}[equation]{Remark}

\newtheorem{hyp}[equation]{Hypotheses}

\newcounter{case}
\newenvironment{case}[1][\unskip]{\refstepcounter{case}\em
\medskip \noindent Case \thecase\ #1.\ }{\unskip\upshape}
\renewcommand{\thecase}{\arabic{case}}

\begin{document}

\title{Strongly Regular Edge-transitive Graphs}

\author{Joy Morris$^1$}\footnote{corresponding author (addresses appear on final page)} 
\address{Department of Mathematics and Computer Science \\
University of Lethbridge \\
Lethbridge, AB. T1K 6R4. Canada}
\email{joy@cs.uleth.ca} 

\author{Cheryl E. Praeger}
\address{School of Mathematics and Statistics \\
University of Western Australia \\
35 Stirling Highway\\
Crawley 6009\\
Western Australia, Australia}
\email{praeger@maths.uwa.edu.au}
\thanks{The second author acknowledges support of a Federation 
Fellowship of the Australian Research Council.}

\author{Pablo Spiga}
\address{Dipartimento di Matematica Pura ed Applicata\\
University of Pad\-ova\\
Via Trieste 63\\
35121 Padova, Italy}
\email{spiga@math.unipd.it}

\begin{abstract}
In this paper, we examine the structure of vertex- and edge-transitive strongly regular graphs, using normal quotient reduction.  We show that the irreducible graphs in this family have quasiprimitive automorphism groups, and prove (using the Classification of Finite Simple Groups) that no graph in this family has a holomorphic simple automorphism group.  We also find some constraints on the parameters of the graphs in this family that reduce to complete graphs. 
\end{abstract}


\maketitle

\section{Introduction}

There has recently been considerable success in using normal quotient 
reduction to analyse the structure of several families of vertex- and 
edge-transitive graphs, including distance-transitive graphs \cite{PSY}, 
$s$-arc-transitive graphs \cite{Praeger-arc-trans,Cheryl}, 
and locally $s$-arc-transitive graphs \cite{GLP-1, GLP-2}.  
This paper initiates a study of strongly regular graphs that are 
vertex- and edge-transitive, using normal quotient analysis.

A \emph{strongly regular graph (srg)} with parameters $(n,k,\lambda,\mu)$ is
a regular graph with $n$ vertices, valency $k$, such that each pair of 
adjacent vertices lies in $\lambda$ triangles, and each pair of non-adjacent
vertices (if such exist) is joined by $\mu$ paths of length 2. In particular
we regard the complete graph $K_n$ as a (somewhat degenerate) 
strongly regular graph 
with parameters $(n,n-1,n-2,0)$.
To simplify the language in this paper, we use the following terminology 
throughout.
\begin{defn}
We refer to a vertex- and edge-transitive strongly regular graph as a {\em ve-srg}.
\end{defn}

Requiring strongly regular graphs to be vertex- and edge-transitive may seem a very stringent condition. However, several important and well-known families of graphs are contained within the class of ve-srgs, making them clearly worthy of study. Notably, Paley graphs, Kneser graphs with $k=2$, and rank 3 graphs are all ve-srgs.  

We begin our study by showing that normal quotient reduction (as defined in Section 2) applies properly to ve-srgs.  That is, we show that each normal quotient of a ve-srg $\Gamma$ is itself a ve-srg.  To ensure that the reduction works properly, we do not require that the group $N$ used in the reduction be normal in the full automorphism group of the graph, but we do require that it be nontrivial, intransitive, and normal in a subgroup of $\Aut(\Gamma)$ that is both vertex- and edge-transitive.  For this reason, a graph that is irreducible must either be complete, or have the property that every vertex- and edge-transitive subgroup of its automorphism group is quasiprimitive (see Definition~\ref{def:qp}).

This reduction sets up our strategy for determining the structure of
ve-srgs: we aim to
\begin{enumerate}
\item show that every connected ve-srg can be reduced to an irreducible
  ve-srg by taking a succession of normal quotients;  
\item characterise and as far as possible determine the irreducible 
ve-srgs: graphs that cannot be  further reduced using normal 
quotient reduction; 
\item use what we learn about the irreducible ve-srgs, and the normal
  quotient structure, to study the structure of
  arbitrary ve-srgs. 
\end{enumerate}
In Section 2 we discuss normal quotient reduction and in particular 
we establish (1). 
In the final three sections of this paper we begin with some minor results 
on (3), and then achieve some more significant progress on (2).

In Section 3, we consider the structure of ve-srgs that reduce to the trivial case of the complete graph under normal quotient reduction.  We completely characterise the structure of ve-srgs that reduce to either $K_2$ or $K_3$, and find 
constraints on the parameters of ve-srgs that reduce to one of the larger 
complete graphs.

In Section 4, we study the family of cartesian products $K_b \square K_b$.  We show that these graphs are ve-srgs, and that graphs in this intriguing family reduce to the complete graph $K_b$ if and only if $b$ is a prime power, 
and are otherwise irreducible.

In Section 5, we show that of the 8 families of quasiprimitive groups, one (holomorphic simple groups) can never arise as the automorphism group of a ve-srg (see Corollary \ref{noHS}).  This result uses the Classification of Finite Simple Groups.

While the word ``connected'' appears in many of the results, it is almost unnecessary; a disconnected strongly regular graph is a disjoint union of cliques of some fixed order, whose structure we completely understand.

\section{Normal quotient reduction}

\noindent
\textbf{Graph Notation:}\quad A graph $\Gamma$ consists of a set $V(\Gamma)$
of vertices
and a subset $E(\Gamma)$ of unordered pairs of vertices called edges. An
automorphism $g$ of $\Gamma$ is a permutation of $V(\Gamma)$ that leaves 
$E(\Gamma)$ invariant. We denote the image of a vertex $x$ under $g$ by $g(x)$.
For a subgroup $H$ of the automorphism group $\Aut(\Gamma)$ of $\Gamma$, we 
denote the orbit of $x$ under $H$ by ${}^Hx=\{h(x)\,|\,h\in H\}$ and we say 
that $\Gamma$ is $H$-vertex-transitive or $H$-edge-transitive if $H$ is 
transitive on $V(\Gamma)$ or $E(\Gamma)$ respectively. We extend this 
usage to possibly unfaithful actions: for example, if $H$ 
acts as a vertex-transitive group of automorphisms of $\Gamma$ 
with kernel $N$, then we will often say that $\Gamma$ is $H$-vertex-transitive,
rather than $(H/N)$-vertex-transitive.

\medskip
In this section, we consider the possibility that $\Gamma$, a
ve-srg,
has a normal quotient.  We look at a transitive subgroup
$G\le \Aut(\Gamma)$, and suppose that there exists some nontrivial
normal subgroup $N$ of $G$ that is intransitive in its action on the
vertices of $\Gamma$. 

\begin{defn}
Let $\Gamma$ be a graph and $G$ a vertex-transitive subgroup
of $\Aut(\Gamma)$.  Suppose that there is some group $N$ such
that $1 \neq N \normal G$, and $N$ is intransitive in its action on
$V(\Gamma)$.   

The {\em quotient graph, $\Gamma_N$,} is the graph whose vertices are
the orbits of $N$, with an edge between two distinct vertices ${}^Nx$
and ${}^Ny$ in $\Gamma_N$, if and only if there is an edge of $\Gamma$
between $x'$ and $y'$, for some $x' \in {{}^Nx}$ and some $y' \in
{{}^Ny}$. 
\end{defn}

Sometimes, the original graph will have a nice covering structure with
respect to the quotient graph. 

\begin{defn}
As before, let $\Gamma$ be a graph, $G\le \Aut(\Gamma)$ a transitive
subgroup of the automorphism group, and $N$ an intransitive normal
subgroup of $G$. 
Suppose that for each pair $\{B, B'\}$ of adjacent $N$-orbits, each
vertex in $B$ is adjacent to exactly $\ell$ vertices in $B'$ (so
$\ell$ is a constant that does not depend on the choice of $\{B,
B'\}$).  Then we say that $\Gamma$ is an {\em $\ell$-multicover} of
$\Gamma_N$. 
\end{defn}

Notice that if $\Gamma$ is an $\ell$-multicover of $\Gamma_N$, then
$\ell$ must divide the valency of $\Gamma$.  We now make some
observations about the structure of the quotient graph. 

\begin{lem}\label{quotient-properties}
Let $\Gamma$ be a connected, vertex- and edge-transitive graph, with
$G\le \Aut(\Gamma)$  acting transitively on the vertices.  Let $1 \neq
N \normal G$, for some intransitive $N$, and let $\Gamma_N$ be the
corresponding quotient graph.  Then  
\begin{enumerate}
\item $\Gamma_N$ is connected, \label{connected}
\item $\Gamma_N$ is vertex-transitive, and \label{vx-trans}
\item $\Gamma_N$ has diameter at most the diameter of $\Gamma$. \label{diam}
\end{enumerate}

Moreover, if $G$ is edge-transitive, then 
\begin{enumerate}\setcounter{enumi}{3}
\item all edges of $\Gamma$ join vertices in distinct
  $N$-orbits,\label{no-edges} 
\item $\Gamma_N$ is $G$-edge-transitive, and\label{edge-trans}
\item $\Gamma$ is an $\ell$-multicover of $\Gamma_N$, for some divisor
  $\ell$ of the valency of $\Gamma$.\label{multicover} 
\end{enumerate}
\end{lem}

\begin{proof}
(\ref{connected} and \ref{diam}) Let ${}^Nx$ and ${}^Ny$ be arbitrary
  vertices of $\Gamma_N$.  Since $\Gamma$ is connected, there is an
  $x$-$y$ path in $\Gamma$, say $x=x_1, x_2, x_3, \ldots, x_i=y$.
  Then for each $j \in \{1, \ldots, i-1\}$, the definition of
  $\Gamma_N$ and the fact that $\{x_j, x_{j+1}\}$ is an edge of $\Gamma$
  tells us that either ${}^Nx_j={{}^Nx_{j+1}}$, or $\{{}^Nx_j, 
{{}^Nx_{j+1}}\}$ is
  an edge of $\Gamma_N$.   Thus if we ignore repeated vertices,
  ${{}^Nx}={}^Nx_1, {{}^Nx_2}, {{}^Nx_3}, \ldots, {{}^Nx_i}={{}^Ny}$ is an
  ${{}^Nx}$\,-\,${{}^Ny}$ path 
  in $\Gamma_N$.  This shows that $\Gamma_N$ is connected.  Moreover,
  this shows that the distance between arbitrary vertices ${}^Nx$ and
  ${}^Ny$ in $\Gamma_N$ is no greater than the distance between $x$ and
  $y$ in $\Gamma$, which implies that the diameter of $\Gamma_N$ is at
  most the diameter of $\Gamma$.  

(\ref{vx-trans}) Let ${}^Nx$ and ${}^Ny$ be arbitrary vertices of
  $\Gamma_N$. Since $G$ acts transitively on the vertices of $\Gamma$,
  there is some $g \in G$ such that $g(x)=y$.  But since $N \normal
  G$, we have $g({}^Nx)={}^Ny$, so $\Gamma_N$ is
  $G$-vertex-transitive. 

(\ref{no-edges}) Since $\Gamma$ is connected and $N$ is intransitive,
  there must be an edge that joins two distinct orbits of $N$.  Since
  $N \normal G$, the orbits of $N$ are blocks of imprimitivity under
  the action of $G$, so no element of $G$ can map this edge to an edge
  both of whose endvertices lie within an orbit of $N$.  But since $G$
  is transitive on the edges of $\Gamma$, this means that there cannot
  be an edge both of whose endvertices lie within any one orbit of
  $N$. 

(\ref{edge-trans}) Let $\{{{}^Nx},{{}^Nx'}\}$ and $\{{{}^Ny},{{}^Ny'}\}$ 
be arbitrary
  edges of $\Gamma_N$. By the definition of $\Gamma_N$, there is some
  $x_1 \in {{}^Nx}$ and some $x_1' \in {{}^Nx'}$ such that $\{x_1,x_1'\}$ is an
  edge of $\Gamma$.  Similarly, there is some $y_1 \in {}^Ny$ and some
  $y_1' \in {{}^Ny'}$ such that $\{y_1,y_1'\}$ is an edge of $\Gamma$. Since
  $G$ acts transitively on the edges of $\Gamma$, there is some $g \in
  G$ such that $g(\{x_1,x_1'\})=\{y_1,y_1'\}$.  But since $N \normal G$,
  we
  have $$g(\{{{}^Nx},{{}^Nx'}\})=g(\{{{}^Nx_1},{{}^Nx_1'}\})=
\{{{}^Ny_1},{{}^Ny_1'}\}=\{{{}^Ny},{{}^Ny'}\},$$   
  so $\Gamma_N$ is $G$-edge-transitive. 

(\ref{multicover}) Let $\{{{}^Nx},{{}^Nx'}\}$ be an arbitrary edge of
  $\Gamma_N$, and suppose that there are edges from $x$ to exactly
  $\ell$ vertices in ${}^Nx'$. Since the action of $N$ is transitive on
  ${}^Nx$ and fixes every $N$-orbit setwise, it follows that every vertex in
  ${}^Nx$ must have edges to exactly $\ell$ vertices in ${}^Nx'$.  So
  there are $\ell |{}^Nx|$ edges between the two $N$-orbits.  Since $G$
  acts transitively on the vertices and $N \normal G$, $|{}^Nx|=|{}^Nx'|$,
  so considering the action of $N$ on ${}^Nx'$ 
shows that every vertex in ${}^Nx'$ has
  edges to exactly $\ell$ vertices in ${}^Nx$.  Now if $\{{{}^Ny}, {{}^Ny'}\}$ 
is
  another edge of $\Gamma_N$, the fact that $\Gamma_N$ is
  $G$-edge-transitive forces every vertex in ${}^Ny$ to be adjacent to
  exactly $\ell$ vertices in ${}^Ny'$. 
\end{proof}

We also deduce that the quotient graph is strongly regular, if the
original graph is. 

\begin{lem}\label{quotient-sr}
Let $\Gamma$ be a connected ve-srg, with $G\le \Aut(\Gamma)$ transitive on the vertices and
the edges.  Let $1 \neq N \normal G$, for some intransitive $N$, and
let $\Gamma_N$ be the corresponding quotient graph.  If $\Gamma$ is
complete, then no such quotient is possible; otherwise, $\Gamma_N$ is
a connected ve-srg. 
\end{lem}

\begin{proof}
We use $(n, k, \lambda, \mu)$ to denote the parameters of $\Gamma$.

That $\Gamma_N$ is connected and vertex- and edge-transitive follows from Lemma~\ref{quotient-properties}(\ref{connected}, \ref{vx-trans} and \ref{edge-trans}).

If $\Gamma$ is complete, then any edge-transitive group $G$ is
2-homogeneous on the vertices of $\Gamma$ (that is, transitive on the
unordered 2-sets).  All transitive 2-homogeneous groups are primitive
(see \cite{DiMo}, page 35), so all of their nontrivial normal
subgroups are transitive. Thus when
$\Gamma$ is complete, there is no nontrivial vertex-intransitive 
group $N$.  We
may therefore assume that $\Gamma$ has diameter 2, as strongly regular
graphs have diameter at most 2. 

If $\Gamma_N$ is complete, then it is strongly regular, and we are done.  So we may
also assume that $\Gamma_N$ has diameter 2 (by
Lemma~\ref{quotient-properties}(\ref{diam}), its diameter is at most
2). 

Since $\Gamma_N$ is edge-transitive, each
of its edges must lie in a constant number of triangles.  This
establishes the parameter $\lambda'$ of $\Gamma_N$. 

Since $\Gamma_N$ is vertex-transitive, it is clear that
$\Gamma_N$ is regular, of degree $k'$ (say). 

Now we wish to show that the number of 2-paths between nonadjacent
vertices of $\Gamma_N$ does not depend on the choice of the
nonadjacent vertices.  Let ${}^Nx$ and ${}^Ny$ be nonadjacent vertices of
$\Gamma_N$.  We count the number of 2-paths in $\Gamma$ between the
sets ${}^Nx$ and ${}^Ny$, in two different ways. 

First we set up some notation for our counting.  Let $b$ denote the
number of vertices of $\Gamma$ in each orbit of $N$.  By
Lemma~\ref{quotient-properties}(\ref{multicover}), $\Gamma$ is an
$\ell$-multicover of $\Gamma_N$ for some divisor $\ell$ of the valency
of $\Gamma$. 

Now, the number of 2-paths in $\Gamma$ between ${}^Nx$ and ${}^Ny$ can be
counted as $b$ choices for a vertex $x_1 \in {{}^Nx}$ to start the path,
times $b$ choices for a nonadjacent vertex $y_1 \in {{}^Ny}$ to end the
path, times $\mu$ 2-paths between $x_1$ and $y_1$, since $\Gamma$ is
strongly regular. 

Alternatively, the number of 2-paths in $\Gamma$ between ${}^Nx$ and
${}^Ny$ can be counted as $b$ choices for a vertex $x_1 \in {{}^Nx}$ to
start the path, times $\mu'_{{}^Nx,{{}^Ny}}$ choices for a set ${}^Nu$ that is
mutually adjacent to ${}^Nx$ and ${}^Ny$, times $\ell$ choices for a
vertex $u_1 \in {{}^Nu}$ that is adjacent to $x_1$, times $\ell$ choices
for a vertex $y_1 \in {{}^Ny}$ that is adjacent to $u_1$. 

These two ways of counting the same thing, show us that 
$$
b^2\mu =b\mu'_{{}^Nx,{}^Ny}\ell^2.
$$  
Since none of $b$, $\mu$ and $\ell$ depends
on the choices of ${}^Nx$ and ${}^Ny$, we have that
$\mu'_{{}^Nx,{}^Ny}=b\mu/\ell^2$ must not depend on the choices of ${}^Nx$
and ${}^Ny$, either.  Thus the number of 2-paths between nonadjacent
vertices of $\Gamma_N$ does not depend on the choice of the
nonadjacent vertices. 
\end{proof}

We have shown that a normal quotient of a connected ve-srg, is itself
a connected ve-srg.

Notice that if $\Gamma$ is a graph that cannot be further reduced
using this normal quotient reduction 
then either $\Gamma$ is complete, or $\Gamma$ has no group $G$ of
automorphisms that is transitive on both the vertices and the edges
and that
has a subgroup $N$ such that $1 \neq N \normal G$ is vertex-intransitive. 
The following definition is therefore very important.

\begin{defn}\label{def:qp}
A transitive permutation group is said to be {\em quasi\-prim\-it\-ive} if 
every nontrivial normal subgroup is transitive. 
\end{defn}
Thus, the ``irreducible" connected ve-srgs are the complete
graphs, together with the graphs for which every vertex- and
edge-transitive group of automorphisms is quasiprimitive.  

It might be thought that the requirement of edge-transitivity is
artificial; it certainly creates some noteworthy awkwardness in at
least two ways.  Firstly, the complement of a connected ve-srg will be strongly
regular, vertex-transitive and often connected, but edge-transitivity
is unlikely.  Secondly, we will see later with the example of the
direct product $K_b \square K_b$ that there may be ve-srgs with normal 
quotients that are also ve-srgs, but even though both
the graph and the quotient are edge-transitive, the reduction can only
be made by taking a normal subgroup of a vertex-transitive subgroup of the
automorphism group, \emph{not} by taking a normal subgroup of a vertex-\ and 
\emph{edge-}transitive subgroup of the automorphism group.  However,
it seems to be very difficult to determine much information about
these graphs if we drop the requirement that $G$ be edge-transitive.   

We will consider first the degenerate case of graphs for which some 
normal quotient is a complete graph, and therefore irreducible.

\section{Complete graphs as normal quotients}

In this section, we consider the degenerate case in which some normal
quotient of a ve-srg is a complete graph.  We analyse the structure of the
original graph in this case. 

We will be using the same hypotheses repeatedly in this section, so to
shorten the statements of the results, we state these hypotheses
here. 

\begin{hyp}\label{hyps}
Let $\Gamma$ be a connected ve-srg with parameters $(n,k,$ $\lambda, \mu)$, and $G\le \Aut(\Gamma)$ acting
transitively on the vertices and on the edges.  Let $1 \neq N \normal
G$, for some intransitive $N$, and let $\Gamma_N$ be the corresponding
quotient graph.  Furthermore, let  $b$ denote the number of vertices
in each orbit of $N$; $m$  the valency
of $\Gamma_N$; and $\ell=k/m$  the number of edges from $x$ to ${}^Ny$
whenever ${}^Nx$ and ${}^Ny$ are adjacent in $\Gamma_N$. 
\end{hyp}

\begin{lem}\label{K2}
Under Hypotheses \ref{hyps}, if $\Gamma_N=K_2$, then $\Gamma$ is a
complete bipartite graph. 
\end{lem}

\begin{proof}
 Let ${}^Nx$ and ${}^Ny$ be the two vertices of $\Gamma_N$.  By
 Lemma~\ref{quotient-properties}(\ref{no-edges}), every edge of
 $\Gamma$ has one endvertex in ${}^Nx$ and the other in ${}^Ny$, so
 $\Gamma$ is bipartite. 

Since $\Gamma$ is strongly regular, it has diameter at most 2. Since
two vertices in opposite sets of the bipartition must be at an odd
distance, any two such vertices must be at distance 1.  Thus, $\Gamma$
is a complete bipartite graph. 
\end{proof}

This result does not generalise, but we can draw some conclusions
about the parameters of ve-srgs with a complete normal quotient. 

\begin{prop}\label{complete-parameters}
Under Hypotheses \ref{hyps}, if $\Gamma_N=K_{m+1}$ for some $m \ge 1$,
then 
\begin{enumerate}
\item $(b-1)\mu = \ell m (\ell-1)$ and in particular $\ell \ge 2$,
\item $(b-\ell)\mu= \ell (\ell m -\ell - \lambda)$,
\item $\mu(\ell -1) = \ell (\ell+\lambda-m),$ and
\item $\ell \mid \mu$,
\end{enumerate}
Furthermore, if $\Gamma$ is not the complete multipartite graph 
$K_{(m+1)[b]}$, then 
\begin{enumerate}\setcounter{enumi}{4}
\item $\mu \le (m-1)\ell$.
\end{enumerate}
\end{prop}

\begin{proof}
By Lemma~\ref{K2}, the result holds for $m=1$.  So in what follows, we
may assume $m \ge 2$. 
By Lemma~\ref{quotient-properties}(\ref{no-edges}), every edge of
$\Gamma$ lies between orbits, so $\Gamma$ is multipartite. 

We count the number of 2-paths that start at a fixed vertex $v$, and
end in ${}^Nv$, in two ways.  First, with $v$ fixed, there are $b-1$
ways to choose a vertex $v' \in {{}^Nv}$, to be the terminal vertex of the
2-path.  Since $v$ and $v'$ are in the same $N$-orbit, they are
nonadjacent by Lemma~\ref{quotient-properties}(\ref{no-edges}), so 
there are $\mu$ different 2-paths from $v$ to $v'$ to
choose amongst.  This makes $(b-1)\mu$ 2-paths in all.  Alternatively,
again with $v$ fixed, we can choose an adjacent vertex $w$ in $k=\ell
m$ ways.  For each such choice, we can choose any of the $\ell-1$
vertices of ${}^Nv$ that are not $v$ but that are adjacent to $w$, to
complete the 2-path.  Thus, we conclude that  
\begin{equation}\label{count1}
(b-1)\mu = \ell m (\ell-1).
\end{equation}
This is the first of our desired conclusions.  Furthermore, since $N
\neq 1$, we have $b>1$, and by Lemma~\ref{quotient-sr}, $\Gamma$ is
not complete, so $\Gamma$ has diameter 2, meaning $\mu \neq 0$.  Thus
Equation \ref{count1} forces $\ell \ge 2$. 

With a fixed starting vertex $v$, we count the number of 2-paths that
start at $v$ and end at some vertex $u \not\in {{}^Nv}$ for which $v$ and
$u$ are not adjacent (if there is no such vertex $u$, our count will produce 0).  Fixing $v$, there are $m$ orbits of $N$ that do
not contain $v$, and since $\Gamma_N=K_{m+1}$, each of these contains
$b-\ell$ vertices that are not adjacent to $v$ that serve as choices
for the terminal vertex $u$ of our 2-path.  For each choice of $u$,
there are $\mu$ 2-paths from $v$ to $u$ from which we may choose.
Thus, the number of 2-paths is $m(b-\ell)\mu$.  Alternatively, again
with $v$ fixed, we can choose an adjacent vertex $w$ in $k=\ell m$
ways.  For each such choice, there are $\ell m- \ell$ vertices that
are adjacent to $w$ but not in ${}^Nv$; however, $\lambda$ of these
vertices are actually adjacent to $v$ also (recall that none of the
vertices in ${}^Nv$ are adjacent to $v$ by
Lemma~\ref{quotient-properties}(\ref{no-edges})), so in total 
$\ell m(\ell m - \ell - \lambda)$ is the number of 2-paths starting at $v$ and
terminating in some $u \not\in {}^Nv$ that is not adjacent to $v$.
We conclude that  
$$
m(b-\ell)\mu= \ell m (\ell m -\ell - \lambda),
$$ 
and dividing through by $m$,
\begin{equation}\label{count2}
(b-\ell)\mu= \ell (\ell m -\ell - \lambda),
\end{equation}
our second conclusion.

Now we take Equation \ref{count1} and subtract Equation \ref{count2}, to obtain
$$\mu(\ell -1) = \ell (\ell+\lambda-m),$$  our third conclusion.
Furthermore, since $\ell$ and $\ell-1$ are coprime, this forces $\ell
\mid \mu$, our fourth conclusion. 

Finally, suppose that $\Gamma$ is not $K_{(m+1)[b]}$, and fix two nonadjacent 
vertices $v,u$ in distinct blocks (such vertices exist since $\Gamma$ is not
complete multipartite). The first vertex $v$ has $(m-1)\ell$ neighbours in
blocks other than ${}^Nu$.  
Since any 2-path between the vertices must go through a
third block (by Lemma~\ref{quotient-properties}(\ref{no-edges})),
the number $\mu$ of 2-paths between $v$ and $u$ cannot be more
than $(m-1)\ell$.  This is the last of our desired conclusions. 
\end{proof}

In the case where $\Gamma$ is not complete multipartite, then, we have
(among other things) that $\ell \mid \mu$, but $\mu \le (m-1) \ell$.
It would seem to be interesting to determine how these relationships
between $\mu$ and $\ell$ affect the graph $\Gamma$, and for which
values of $\mu$ relative to $\ell$ there in fact are strongly regular,
connected, edge-transitive graphs.  We begin by considering one
extreme of the relation: the case where $\mu = (m-1)\ell$. 

\begin{prop}\label{mu=m-1l}
Under Hypotheses \ref{hyps}, if $\Gamma_N=K_{m+1}$ for some $m\geq2$ and
$\mu=(m-1)\ell$, then $\Gamma \cong K_{b[b]}-bK_b$, the complete
$b$-partite graph with $b$ vertices in each part, with the edges of
$b$ vertex-disjoint copies of $K_b$ deleted, and $b=\ell +1=m+1$. 
\end{prop}

\begin{proof}
The graph $\Gamma$ cannot be
$K_{(m+1)[b]}$ since for that graph $\mu=m\ell$. Therefore
Proposition~\ref{complete-parameters}(1) gives 
$(b-1)\mu=\ell m(\ell-1)$, so $(b-1)(m-1)=m(\ell-1)$. 
Thus $b-1=(1+\frac{1}{m-1})(\ell-1)\le 2(\ell-1)$, so $b \le 2\ell -1$.

Fix a vertex $v$. If $b>\ell+1$, then there are at least two vertices
$u$ and $w$ in some other block, both of which are nonadjacent to $v$.
Since $\mu=(m-1)\ell$, both $u$ and $w$ must have exactly the same
neighbours as $v$ in the remaining $m-1$ blocks, to create $\mu$
2-paths between $v$ and $u$, and $\mu$ 2-paths between $v$ and $w$.
Since these $\mu$ mutual neighbours also create $\mu$ 2-paths between
$u$ and $w$, it must be that $u$ and $w$ have no mutual neighbours in
the block containing $v$.  But since each has $\ell$ neighbours in the
block containing $v$, and there are a total of $b\le2\ell-1$ vertices
in that block, the pigeon-hole principle forces $u$ and $w$ to have at
least one mutual neighbour in the block containing $v$, a
contradiction.  We can thus conclude that $b \le \ell+1$.  In fact,
since $b \ge \ell$ by definition of $\ell$, and $b=\ell$ would give a
complete multipartite graph, we can conclude that $b=\ell+1$. 

We claim that  any two vertices $v$ and $u$ in the same block have $\ell-1$
common neighbours in each of the other $m$ blocks. Since $\Gamma$ is an 
$\ell$-multicover of $\Gamma_N$, and since $b=\ell+1$, the number of common 
neighbours in each block is at least $\ell-1$.  If $v$ and $u$ had 
$\ell$ common neighbours in some block, then  the other vertex in that 
block could have at most $b-2=\ell-1$
neighbours in ${}^Nv$, which is a contradiction.  Thus the claim is proved, 
and hence $\mu=m(\ell-1)$.  By
assumption, then, $(m-1)\ell = m(\ell-1)$, which forces $m=\ell=b-1$.
We also have $\mu=\ell(\ell-1)$ and $\lambda=(\ell-1)^2$ (from
Proposition~ \ref{complete-parameters}, parts (1) and (3)).  These parameters force
the structure that we have claimed, as we will show in the next
paragraph. 

Towards a contradiction, suppose that for some vertex $v$, the set of
non-neighbours of $v$ that are not in ${}^Nv$ is not an
independent set of vertices.  Then there are two non-neighbours of
$v$, $u$ and $w$ (say), that are adjacent, and $u$, $v$ and $w$ are
all in different blocks.  Now, since each of $u$ and $w$ has
$\ell=b-1$ neighbours in each of the blocks that do not contain $u$ or
$w$, they must have either $\ell$ or $\ell-1$ common neighbours in
each of these blocks.  Since there are $\ell-1$ such blocks, and they
have a total of $(\ell-1)^2$ mutual neighbours, they must have exactly
$\ell-1$ mutual neighbours in each of these blocks.  But since both $u$ and $w$
are non-neighbours of $v$, they have $\ell$ mutual neighbours in the
block ${}^Nv$, a contradiction. 
\end{proof}

It is not hard to verify that the graphs $K_{b[b]}-bK_b$ form an
infinite family of connected ve-srgs
that are not complete multipartite, and that have complete normal 
quotients.

Interestingly, a consequence of the preceding result is that there is 
only one ve-srg that
is not complete tripartite and has $K_3$ as a normal quotient. 

\begin{cor}\label{K_3}
Under Hypotheses \ref{hyps}, if $\Gamma_N=K_{3}$, then either $\Gamma
\cong K_{3[b]}$, or $b=3$ and $\Gamma \cong K_{3[3]}-3K_3$.
\end{cor}

\begin{proof}
Here we have $m=2$.
We assume that $\Gamma$ is not $K_{3[b]}$.

Since $m=2$, Proposition~\ref{complete-parameters}(4) and (5) give
$\ell \mid \mu$ (so $\ell \le \mu$), and $\mu \le \ell$.  Thus
$\mu=(m-1)\ell=\ell$.  Now Proposition~\ref{mu=m-1l} completes the
proof. 
\end{proof}

At the other extreme lies the possibility that $\mu=\ell$.  Here there
is another infinite family of connected ve-srgs: the graphs $K_b \square K_b$.  As there is
quite a bit to say about these graphs, we will discuss them in a
separate section, shortly. 

Meanwhile, putting additional standard facts about strongly regular 
graphs together with the results of Proposition~\ref{complete-parameters}, 
we can produce a bit more information about the possible parameters 
of other ve-srgs $\Gamma$ that reduce to complete graphs. By the eigenvalues of
$\Gamma$ we mean the eigenvalues of its adjacency matrix (see \cite{GoRo}).

\begin{prop}
Under Hypotheses \ref{hyps}, if $\Gamma_N=K_{m+1}$, then either
$\Gamma$ is a complete multipartite graph, or the eigenvalues of
$\Gamma$ are  
$$k\text{, }\theta=m-r\text{, and }\tau=-\ell,$$
where $\ell=k/m$ and $r=\mu/\ell$.  Furthermore, the multiplicities of
the eigenvalues $\theta$ and $\tau$ are 
\begin{eqnarray*}
m_{\theta}&=& \frac{m\ell(m+1)(\ell-1)}{m-r+\ell},\\
m_{\tau}&=&\frac{m}{m-r+\ell}(m^2\ell-m^2+2rm+m\ell-m+\ell rm-r^2+r),
\end{eqnarray*}
so these values must be integers.

The parameters for $\Gamma$ are $$\bigl((m+1)(m\ell-m+r)/r,m\ell,
(r-1)\ell+m-r, r\ell\bigr).$$ 
\end{prop}

\begin{proof}
We assume that $\Gamma$ is not multipartite, and deduce the given
formulae for the eigenvalues, their multiplicities, and the parameters
of $\Gamma$. 
Standard results (cf. \cite{GoRo}, p. 220) on strongly regular graphs
give the formulas 
$$\theta=(\lambda-\mu+\sqrt{\Delta})/2\text{ and
}\tau=(\lambda-\mu-\sqrt{\Delta})/2,$$ where
$\Delta=(\lambda-\mu)^2+4(k-\mu)$. 
Solving for $\lambda$ in Proposition~\ref{complete-parameters}(3) gives 
$\lambda=(r-1)\ell+m-r$
(which is the value given above for the third parameter of $\Gamma$).
Substituting this into the formula for $\Delta$ gives
$\Delta=(m+\ell-r)^2$.  The values given for $\theta$ and $\tau$ are
an immediate consequence of combining these results. 

Similarly, standard results on strongly regular graphs give
$$m_{\theta}=\frac{(n-1)\tau+k}{\tau-\theta}\text{ and
}m_{\tau}=\frac{(n-1)\theta+k}{\theta-\tau}.$$ 
Solving for $b$ in Proposition~\ref{complete-parameters}(1) gives
$b=(m\ell -m+r)/r$, and we know $n=(m+1)b$.  This yields the given value for
$n$ in the parameters of $\Gamma$, completing the calculations of the
parameters.  Using this value for $n$ in the formulas for $m_{\theta}$
and $m_{\tau}$ completes the result. 
\end{proof}

While we have determined some information about ve-srgs that reduce to
the complete graph, and have determined some special cases, there is
much still to be explored here. 

\section{The family $K_b \square K_b$}

In this section, we consider the infinite family of graphs consisting
of cartesian products of two copies of a complete graph.  First we
show that these graphs are all connected ve-srgs.  
Then we show that while all have
$K_b$ as a normal quotient relative to a vertex-transitive subgroup of 
automorphisms, the family is divided into two infinite
subfamilies. The graphs $\Gamma$ in the first subfamily have an
edge-transitive group $G$ of automorphisms and a nontrivial normal
subgroup $N$ of $G$, such that $\Gamma_N= K_b$.  In the other
subfamily, every edge-transitive group of automorphisms is
quasiprimitive, so graphs in the second subfamily are irreducible. 

We define $\Gamma=K_b\square K_b$ as the graph with vertices the 
pairs $(i,j)$ with $i,j \in \Z_b$, such that
$(i,j)$ and $(i',j')$ form an edge if and only if $i=i'$ or 
$j=j'$ (but not both).

\begin{prop}\label{infinite-family}
Let $\Gamma = K_b \square K_b$.  Then
$\Gamma$ is vertex- and edge-transitive and strongly regular with
parameters $(b^2, 2b-2, b-2, 2)$. 

Furthermore, 
$\Gamma$ is a $b$-partite
graph with parts of cardinality $b$.
In fact, $K_b$ is a normal quotient of $K_b \square
K_b$ relative to a vertex-transitive subgroup of $\Aut(\Gamma)$, and
$\Gamma$ is a $2$-multicover of this quotient. 
\end{prop}

\begin{proof}
It is well-known and easy to show that $\Gamma$ is vertex- and edge-transitive.
The number of vertices is $b^2$ and the valency is $2b-2$.
Every edge is in a unique $K_b$, either formed by the $b$ vertices
with the same first coordinate, or by the $b$ vertices with the same
second coordinate. No edge is in any other triangles. This establishes
that $\lambda=b-2$. 

If $(i,j)$ and $(i', j')$ are not adjacent, then $j \neq j'$ and $i
\neq i'$.  The only vertices adjacent to both
of these vertices, are $(i,j')$ and $(i',j)$. This establishes that
$\mu=2$ and that $\Gamma$ is indeed strongly regular. 

Choose the diagonal sets $\bigl\{\{(i,i+j): i \in \Z_b\}:j \in
\Z_b\bigr\}$ as the partition sets. Then it is clear that $\Gamma$ is
$b$-partite with parts of cardinality $b$, and each vertex is joined 
to exactly two vertices in each partition set apart from the one containing it.

Let $G$ be $\Z_b \times \Z_b$ in its action by left
multiplication on the vertices of $\Gamma$. Then $G$ is a subgroup
of $\Aut(\Gamma)$ that acts transitively on the vertices.
Furthermore, as $G$ is abelian, the subgroup $N=\{n_x: x\in \Z_b\}$,
where $n_x(i,j)=(i+x,j+x)$, is normal in $G$, and the
$N$-orbits are the partition sets mentioned above.  Thus $\Gamma_N\cong K_b$
and $\Gamma$ is a $2$-multicover of $\Gamma_N$. 

We can also obtain $K_b$ as a normal quotient by taking $G$ to be $S_b \times
S_b$, which again is vertex-transitive but not edge-transitive, and $N
= S_b \times 1$. However in this case $\Gamma$ is not a multicover of 
$\Gamma_N$, since the partition sets this time contain edges of $\Gamma$. 
\end{proof}

Unfortunately, an unexpected artifact of our requirement that the
subgroup $G$ of $\Aut(\Gamma)$ used in these results must be
transitive on the edges as well as on the vertices, is that some, but
not all, of the graphs in this family actually reduce to the complete
graph using our normal quotient reduction scheme, even though both the
original graphs in this family and the quotients are edge-transitive. 

\begin{prop}\label{K_bcartesian}
Let $\Gamma$ be the cartesian product $K_b \square K_b$.  If $b$ is a
prime power, then there is a vertex- and edge-transitive sugroup 
$G \le \Aut(\Gamma)$ with a normal
subgroup $N$ such that $\Gamma_N\cong K_b$. 
\end{prop}

\begin{proof}
Let ${\rm GF}(b)$ denote the field of order $b$, and label the vertices 
of $\Gamma$ by pairs $(i,j)$ with $i,j \in {\rm  GF}(b)$ such that 
$(i,j)$ and $(i',j')$ form an edge if and only if $i=i'$ or $j=j'$ (but
not both). 

Let $n_x(i,j)=(i+x,j+x)$ for every $i, j, x \in {\rm GF}(b)$, and let
$N=\{n_x: x \in {\rm GF}(b)\}$.  Then $N \le \Aut(\Gamma)$.  The
orbits of $N$ are the $b$ transversals $\{(i+x, x): x \in {\rm
  GF}(b)\}$, one for each value of $i$. Also, as in the proof of
Proposition~\ref{infinite-family}, $N$ is a normal subgroup of a
vertex-transitive subgroup of $\Aut(\Gamma)$, and $\Gamma_N=K_b$.  It
remains to be shown that the normaliser of $N$ in $\Aut(\Gamma)$ is in
fact edge-transitive. 

Let $g_{r,s,s'}(i,j)=(ri+s, rj+s')$ for every $i, j, s, s' \in {\rm
  GF}(b)$, and for every $r \in {\rm GF}(b)^*$.  Let $G' =\{g_{r, s,
  s'}:  s, s' \in {\rm GF}(b), r \in {\rm GF}(b)^*\}$.  Let
$\delta(i,j)=(j,i)$ for every $i, j \in {\rm GF}(b)$, and let
$G=\langle G', \delta \rangle$.  Then it is straightforward to see
that $N \normal G$ and $G\leq \Aut(\Gamma)$.  

We claim that $G$ acts transitively on the vertices and edges of
$\Gamma$.  By varying $s$ and $s'$ with $r=1$, it is clear that $G$
acts transitively on the vertices of $\Gamma$.  The group generated by
the automorphisms $g_{r,0,0}$ together with $\delta$ fixes the vertex
$(0,0)$ and acts transitively on its neighbours, so $G$ does indeed
act transitively on the edges of $\Gamma$. 
\end{proof}

Our next series of results shows that in fact, the graphs $K_b \square
K_b$ reduce to $K_b$ only if $b$ is a prime power, and that otherwise,
these graphs are irreducible. 

\begin{lem}\label{Norbits}
Let $\Gamma$ be the cartesian product $K_b \square K_b$.  If $G\le
\Aut(\Gamma)$ acts transitively on the vertices and the edges, then
the orbits of any nontrivial intransitive normal subgroup of $G$ must: 
\begin{itemize}
\item be independent sets in $K_b \square K_b$; 
\item have length $b$; and
\item yield a complete normal quotient $\Gamma_N\cong K_b$. 
\end{itemize}
\end{lem}

\begin{proof}
Let $N$ be a nontrivial intransitive normal subgroup of $G$.  Then the
$N$-orbits in $V(\Gamma)$ form blocks of imprimitivity for the action of $G$.
Since $G$ is edge-transitive and $\Gamma$ is connected, there cannot
be any edges between vertices that lie in the same (nontrivial) block
of imprimitivity of $G$.  This establishes that since $N$ is
intransitive, the orbits of $N$ are independent sets (partial
transversals) of $K_b \square K_b$. 

Since $N$ is nontrivial, its orbits must have length greater than 1.
Let $(x,y)$ and $(x',y')$ be distinct vertices in one orbit of $N$.
These two vertices are not adjacent and so $x\ne x', y\ne y'$.
Since this orbit is an independent set, $(x,y')$ is in a different
$N$-orbit, and these orbits are clearly adjacent in $\Gamma_N$.
Furthermore, since ${}^N(x,y)$ is independent, $(x,y)$ and $(x',y')$ are 
the only vertices of ${}^N(x,y)$ with first entry $x$ or second entry $y'$, 
respectively, and hence they are the only vertices of ${}^N(x,y)$ which
are adjacent to $(x,y')$.  Since $\Gamma$ is an
$\ell$-multicover of $\Gamma_N$, this means that $\ell=2$.  Now, if
${}^N(x,y)$ has length less than $b$, then there is some $j\in{\rm GF}(b)$ 
such that no vertex of ${}^N(x,y)$ has second entry $j$.  Hence $(x,j)$ is 
not in ${}^N(x,y)$, and has exactly one neighbour $(x,y)$ in ${}^N(x,y)$, 
contradicting $\ell=2$.  Thus, the orbits ${}^N(x,y)$ must have length $b$, 
and must contain exactly one pair from each row and each column. It follows 
that the normal quotient $\Gamma_N$ is a complete graph $K_b$.
\end{proof}

\begin{prop}\label{bfactors}
Let $\Gamma$ be the cartesian product $K_b \square K_b$.  Let $G \le
\Aut(\Gamma)$ act transitively on the vertices and the edges of
$\Gamma$.  If $b$ is not a prime power, then $G$ is quasiprimitive. 
\end{prop}

\begin{proof}
The automorphism group of $K_b \square K_b$ is $S_b \wr S_2$, and is
generated by $B=S_b \times S_b$, together with $\delta$, where
$\delta((i,j))=(j,i)$.  Notice that $B \normal \Aut(\Gamma)$; in fact,
$B$ has index 2 in $\Aut(\Gamma)$.  Furthermore, each row and column of
$\Gamma$ is a block of imprimitivity for $B$, so for any subgroup of
$B$, it is reasonable to talk about the action on the rows and
on the  columns of $\Gamma$. 

Let $G$ be a subgroup of $\Aut(\Gamma)$ that acts transitively on the
vertices and edges of $\Gamma$, and let $G'= G \cap B$.  Since $G$ is
edge-transitive, $G$ cannot be contained in $B$. So $G'$ has index $2$
in $G$. Furthermore, the action of $G'$ is 2-homogeneous on the rows
and on the columns of $K_b \square K_b$. 

Let $N \normal G$, with $N$ nontrivial and intransitive.  Then $N' =
N\cap B \normal G'$, and $N'$ has index 1 or 2 in $N$.  By
Lemma~\ref{Norbits}, the orbits of $N$ on the rows must have length
$b$, so the orbits of $N'$ have length $b$ or $b/2$.  Since $G'$ acts
2-homogeneously (and hence primitively) on rows and columns, the
orbits of $N'$ on the rows must have length 1 or $b$. 
Since $b$ is not a prime power, $b \neq 2$, so the orbits of $N'$ on
the rows have length $b$; that is, $N'$ is transitive on the rows.
Similarly, $N'$ is transitive on the columns. 

Suppose that for some $n \in N'$, the action of $n$ fixes row $i$.
Then for any column $j$, $n((i,j))=(i,j)$ since the orbits of $N'$,
being subsets of the orbits of $N$, are independent sets (by
Lemma~\ref{Norbits}).  Thus $n$ fixes every column setwise.
Similarly, since $n$ now fixes column $j$, $n$ fixes every row
setwise.  This shows that the action of $N'$ is faithful and 
regular on the rows.
Then, since $N' \normal G'$  
and $G'$ is 2-homogeneous on the rows, it follows that $N'$ is elementary 
abelian, and in particular $b$ is a prime power.


This shows that $G$ cannot have a nontrivial, intransitive normal
subgroup if $b$ is not a prime power.  By definition of
quasiprimitivity, $G$ must be quasiprimitive if $b$ is not a prime
power. 
\end{proof}

Combining the preceding results yields the following information about
this family of graphs. 

\begin{cor}
Let $\Gamma$ be the cartesian product $K_b \square K_b$. There exists
a vertex- and edge-transitive group $G\le\Aut(\Gamma)$ with a 
nontrivial, vertex-intransitive, normal subgroup $N$ if and only if $b$ 
is a prime power. Moreover in this case each such subgroup $N$
corresponds to a complete normal quotient $\Gamma_N\cong K_b$.
\end{cor}

\begin{proof}
If $b$ is a prime power, the existence of the groups $G$ and $N$ follows from
Proposition~\ref{K_bcartesian}, and the fact that each intransitive $N$
corresponds to a complete normal quotient
$\Gamma_N\cong K_b$ follows from Lemma~\ref{Norbits}.  If $b$ is
not a prime power, the quasiprimitivity of $G$ from
Proposition~\ref{bfactors} shows that no such $N$ exists. 
\end{proof}

Juxtaposing the fact that $K_b$ is always a normal quotient of $K_b
\square K_b$ relative to a vertex-transitive group $G$, with the fact
that $G$ can be edge-transitive if and only if $b$ is a prime power, 
illustrates the complex behaviour of these Cartesian products.
%
Our results show that either these graphs are themselves irreducible,
or they reduce to a complete graph $K_b$ (which is irreducible) under 
normal quotient reduction. 

\section{Irreducible graphs and Holomorphic Simple groups}

By our reduction method, if a graph that is not complete is to be irreducible, then every
group of automorphisms that acts vertex- and edge-transitively on the graph must
be quasiprimitive. 

According to Praeger's characterisation of quasiprimitive groups
\cite{Cheryl}, they fall into one of the following families: 
\begin{enumerate}
\item holomorphic simple groups;
\item affine groups;
\item almost simple groups;
\item simple diagonal action groups with one minimal normal subgroup;
\item compound holomorphic groups;
\item compound diagonal action groups;
\item product action groups; and
\item twisted wreath action groups.
\end{enumerate}

While some of these groups do act as automorphism groups of ve-srgs 
(such as affine groups acting on the Paley graphs, and certain simple group
actions on ve-srgs), we
will devote the remainder of this paper to proving that holomorphic
simple groups cannot arise as vertex- and edge-transitive 
automorphism groups of ve-srgs. This result is stated
in Corollary \ref{noHS}.

We recall the definition of homomorphic simple groups. Let $T$ be a
nonabelian simple group and $W=T\rtimes \Aut (T)$. The group $W$ has a
natural action on the set $T$, where the subgroup $T$ of $W$
acts by left multiplication and the subgroup $\Aut (T) $ of $W$ acts by
automorphisms. With this permutation representation, $W$ is a primitive
subgroup of $\Sym(T)$. We note that the stabilizer of the point $1_T$
in $W$  is $\Aut(T)$. Any subgroup $G$ of $W$ containing $T\rtimes 
\Inn T$ is said to be a holomorphic simple group. In this case, the
group $G$ is quasiprimitive (it is actually primitive), with exactly
two minimal normal subgroups, $T$ and $M=C_G(T)$. Furthermore, $T\cong
M$ and $T\rtimes \Inn T=T\times M$.

\begin{hey}\label{Cayley}
If a holomorphic simple group $G$ (as above) has a vertex-transitive
action on a graph $\Gamma$, then (by definition of quasiprimitive) $T$
and $M$ act transitively on the vertices of $\Gamma$.  In fact, $T$
and $M$ act regularly on the vertices of $\Gamma$  (the subgroups $T$
and $M$ of $G$ are the left and right regular permutation
representations of the simple group $T$).  Thus
$\Gamma$ is a Cayley graph on $T$ (or on $M$). 
\end{hey}

This remark shows that if $G$ is a holomorphic simple
group acting vertex-transitively on $\Gamma$, we can identify the
vertices of $\Gamma$ with the elements of the simple group $T$.   

Again, we have a series of hypotheses that we will be using in many of
the results in this section, so we collect them here.  

\begin{hyp}\label{HS-hyps}
Let $\Gamma$ be a connected ve-srg with parameters $(n,k,$ $\lambda,\mu)$, 
and with $G\le \Aut(\Gamma)$  acting transitively
on the vertices and on the edges.  Suppose that $G$ is a holomorphic
simple group, with minimal normal subgroups $T$ and $M$.  Identify the
vertices of $\Gamma$  with the elements of $T$, according to
Remark~\ref{Cayley}.  Let $H=G_{1_T}$, the stabiliser of the identity, and
let  $S$ be the connection set of the Cayley graph $\Gamma$.   Let $y=\gcd(|H:\Inn T|, |T|-1)$.
\end{hyp}


We need to consider what the connection set of such a Cayley graph
might look like.  The following result is a special case of
Proposition~1 of \cite{Cheryl-normal-Cayley}, but as the proof is
short we include it here. 

\begin{lem}\label{connection}
Under Hypotheses~\ref{HS-hyps}, either $S$ is an orbit of $H$, or for
every $s \in S$ there is no $h \in H$ such that $h(s)=s^{-1}$.  In the
latter case, $S$ is the union of two orbits $B$ and $B^{-1}$ of $H$. 
\end{lem}

\begin{proof}
Suppose $s, s' \in S$ so that $\{1_T,s\}$ and $\{1_T,s'\}$ are edges 
of $\Gamma$.  Since $G$ is edge-transitive, there is either
some $h \in H$ such that $h(s)=s'$, or some $g\in G$ such that $g$
maps $s$ to $1_T$ and $1_T$ to $s'$.  But as $G=HT$, if $g=ht$ with
$t \in T$ and $h \in H$, we must have $t=s^{-1}$ and $h$ takes
$s^{-1}$ to $s'$.  Thus, $s'$ is in the same $H$-orbit as either $s$
or $s^{-1}$. 
\end{proof}

Now we can find an upper bound on the order of the centraliser of any
element in the connection set. 

\begin{lem}\label{centraliser}
Under Hypotheses~\ref{HS-hyps}, we have $\Gamma\not\cong K_n$, 
$y \ge 5$, $y$ is odd and not divisible by 3, and 
\begin{enumerate}
\item if $S$ consists of a single $H$-orbit, then for any element $s
  \in S$, \label{1orbit} 
$$|C_H(s)| \le (y+1)|H:\Inn T|;$$
\item if $S$ is the union of two $H$-orbits, then for any element $s
  \in S$, \label{2orbits} 
$$|C_H(s)| \le 2(y+1)|H:\Inn T|.$$ 
\end{enumerate}
\end{lem}

\begin{proof}
Suppose that $\Gamma\cong K_n$. Since $G$ is edge-transitive, $G$ must 
be 2-homogeneous on $V(\Gamma)$ (that is, transitive on unordered pairs 
of vertices). However, all finite 2-hogeneous groups are primitive of affine
or almost simple type (see \cite[Theorems 4.1B and 9.4B]{DiMo}). 
Since $G$ is holomorphic simple we conclude that $\Gamma\not\cong K_n$.

By Lemma~\ref{connection}, the elements of $S$ form either one orbit
$B=B^{-1}$ (since $S$ is inverse-closed) of $H$, or two
orbits $B, B^{-1}$.  Let $k=|S|$.  Then if we fix $s \in S$, we
have either $S=\{h(s): h \in H\}$ or $S=\{h(s): h \in H\} \cup
\{h(s^{-1}) : h \in H\}$, and $h(s)=s$ if and only if $h \in C_H(s)$.
Thus for any $s \in S$, we have $k=|H:C_H(s)|$ if $S$ is a single
$H$-orbit, or $k=2|H:C_H(s)|$ if $S$ consists of two $H$-orbits. 

Now, let $d$ be the number of orbits of $\Inn T$ on $B$. Since 
$\Inn T $ is normal in $H$, and $H$ is transitive on $B$, it follows that
$H$ permutes the $\Inn T$-orbits in $B$ transitively. In particular,
\begin{equation}\label{dHT}
d \mid |H:\Inn T|
\end{equation}
and all $\Inn T$-orbits on
$S$ have the same length; let $k'$ be that length.
 Then for any $s \in S$, we have $k'=|\Inn T:C_{\Inn T}(s)|$, and
 $dk'=|B|=|H:C_H(s)|$. 

We count the number of edges $(s,s')$, where $s \in S$ and $s'\not\in
S$ with $s' \neq 1_T$, in two different ways.  First, there are $k$
choices for $s\in S$; by the structure of a strongly regular graph,
each of these has $k-\lambda-1$ neighbours outside of $S \cup \{1_T\}$.
Second, there are $|T|-k-1$ vertices outside of $S \cup \{1_T\}$ from
which to choose $s'$; by the structure of a strongly regular graph,
each of these is at distance 2 from the vertex $1_T$, and there are
$\mu$ common neighbours of $s'$ and $1_T$ from which to choose $s$.
Therefore, 
\begin{equation}\label{kbmu}
k(k-\lambda -1)=\mu(|T|-k-1).
\end{equation}
This means either $$dk'(k-\lambda-1)=\mu(|T|-k-1)$$
or $$2dk'(k-\lambda-1)=\mu(|T|-k-1)$$ (depending on whether $S$
consists of one or two $H$-orbits). As $k'\mid |T|$ and $k' \mid k$,
we have $\gcd(k', |T|-k-1)=1$, so $\gcd(k, |T|-k-1)=\gcd(k,|T|-1)$
must divide both $k/k' \in \{d, 
2d\}$, and $|T|-1$.  In fact, all nonabelian simple groups have even
order, so $\gcd(k, |T|-k-1)$ always divides $d$ and $|T|-1$.
Hence,~\eqref{kbmu} yields
\begin{equation}\label{tk1dividesb}(|T|-k-1) \mid \gcd(d,
  |T|-1)(k-\lambda-1).
\end{equation} 
In particular, $|T|-k-1 \le \gcd(d, |T|-1) (k-\lambda-1)$.

Now, $k-\lambda-1 \le k-1$,
so $$|T|\le (\gcd(d, |T|-1)+1)k-\gcd(d, |T|-1)+1\le(\gcd(d,
|T|-1)+1)k.$$  Hence, $$k\ge |T|/(\gcd(d, |T|-1)+1).$$ If   
$k=|H:C_H(s)|$, that is, if $S$ is an $H$-orbit, we see that 
$$
|C_H(s)| \le |H:\Inn T|(\gcd(d, |T|-1)+1),
$$
and part (1) is proved.  If $k=2|H:C_H(s)|$,
then similarly we obtain part (2).

Finally, we show that $y \ge 5$ and that $y$ is not divisible by 3.  

In the special case that $y=1$ we have $\gcd(d,|T|-1)=1$
so~\eqref{tk1dividesb} yields $(|T|-k-1) \mid
(k-\lambda-1)$.  
Since $\mu \le k$ (the number of 2-paths between two vertices cannot
be greater than the valency of each vertex),~\eqref{kbmu}
forces $k=\mu$.  But this would mean that each vertex at distance 2
from the vertex $1_T$ has exactly the same neighbours as the vertex
$1_T$. Now, it is not hard to see that the set of all vertices that have exactly the same neighbours as $1_T$, forms a block
of imprimitivity of $\Aut(\Gamma)$. Therefore, as $G$ is primitive,
we get that there are no vertices at distance 2 from the vertex $1_T$, in
which case $\Gamma$ is complete, which is a contradiction.  So $y>1$.   

Since $|T|$ is even for every nonabelian simple group, $y$ is odd and
in particular $y \neq 2, 4$. 
Suppose that $y=3$. Then $|T|$ is coprime to $3$. The
only nonabelian simple groups whose orders are not divisible by 3 are
the Suzuki groups. So, $T={^2B_2}(q)$
for some $q=2^{2a+1}$. Now $$|T|=q^2(q^2+1)(q-1) \equiv 1\cdot 2\cdot
1 \equiv 2 \pmod{3}.$$ 
Thus $3 \nmid |T|-1$ for any nonabelian simple group $T$, so $y\neq 3$.  Thus $y
\ge 5$, and $3 \nmid y$. 
\end{proof}

In the next result, we quickly dispose of the possibility that the nonabelian simple group $T$ is isomorphic to an alternating or sporadic simple group.  The remainder of the paper will be devoted to eliminating the possibility that $T$ is isomorphic to a finite simple group of Lie type. 
We include the Tits group
$({^2F_4}(2))'$ in the list of sporadic groups, as $({^2F_4(2)})'$ is not properly a simple group of
Lie type.

\begin{thm}\label{alternating}
Under Hypotheses~\ref{HS-hyps},  the group $T$ is not an
alternating or a sporadic group.  
\end{thm}

\begin{proof}
For every alternating or sporadic group $T$, we have $|\Out T|
\in \{1,2,4\}$.  Thus we must have $|H:\Inn T| \in 
\{1,2,4\}$.  The order of every simple group is even, so
$y=1$.  But this contradicts Lemma~\ref{centraliser}. 
\end{proof}

Now we begin to deal with finite simple groups of Lie type.
Our first step involves presenting a general lower bound on the order of the
centralisers in all of these groups.   

\begin{hey}
To state the following lemma so as to include the Ree and Suzuki groups,
the notation that 
we use for some of the finite simple groups of Lie type is not the
most common.  As the notation is not completely standard, we
feel at liberty to do this. We denote a simple group of Lie type by
${^rL_n}(t)$, where $L=A,\ldots,G$ is the Lie type, $n$ is the rank,
$r$ is the order of the graph automorphism of the corresponding
Dynkin diagram and $t$ is the size of the field where the group is
defined. For instance, $^2A_2(q^2)$ is the unitary group of order
$q^3(q^2-1)(q^3+1)/\gcd(3,q+1)$.  
\end{hey}

\begin{lem}\label{Lie-lower-bd}
For any simple group $T= \ ^rL_n(t)$ of Lie type and any element $x$
of the group $T$, we have $|C_T(x)|\ge (t^{1/r}-1)^n/d$, where $d$ is
the order of the diagonal multiplier (see~\cite[page xvi, Table~$6$]{Atlas}).  
\end{lem}

\begin{proof}
Let $G$ be the connected simple algebraic group of rank $n$, of Lie type $L$, of
adjont isogeny type and over the algebraic closure of the field
$\mathbb{F}_t$ of order $t$. Let
$\sigma$ be the Lang-Steinberg endomorphism of $G$ with fixed point group
$G_\sigma$ such that $T\subseteq G_\sigma$. Let $x$ be an element of $T$. By
Lemma~$3.4$ in~\cite{FuGu}, we have that $|C_{G_\sigma}(x)|\geq (t^{1/r}-1)^n$.
Now, since $T$ has index $d$ in $G_\sigma$, we have
that $|C_T(x)|\geq |C_{G_\sigma}(x)|/d\geq (t^{1/r}-1)/d$. This
completes the proof. 
\end{proof}

Before proving the main result of this section, we prove a lemma about
simple groups of Lie type that we will use in our main
proof. 

\begin{lem}\label{PSL-PSU-gcd}
Let $T={^rL_n}(t)$, where $t=p^a$ for some prime $p$, and
$\hat y=\gcd(|\Out T|, |T|-1)$. If $T$
is not a Ree or a Suzuki group, then $\hat y\mid (a/r)$. If $T$ is a Ree or a
Suzuki group, then $\hat y\mid a$.  
\end{lem}

\begin{proof}The order of the automorphism group of $T$ can be found
  in~\cite{Atlas} page xv and is tabulated in Table~$5$. In
  particular $|\Out T|=dfg$, where $d$ is the
  diagonal multiplier, 
  $f$ is the order of the field automorphisms and $g$ is the order of the
  graph automorphisms of the Dynkin diagram (modulo field
  automorphisms). By inspection, $d,g$ 
  divide the order of $T$ and $f=a$. Also, by inspection of Table~$5$
  in~\cite{Atlas}, we 
  have that $r$ divides $|T|$ and $a$ if $T$ is not a Ree or a Suzuki
  group. Therefore the lemma follows.   
\end{proof}

\begin{thm}\label{Lie-type}
Under Hypotheses~\ref{HS-hyps}, the group $T$ is not a finite simple
group of Lie type.
\end{thm}

\begin{proof}
Assume $T={^rL_n}(t)$, where $t=p^a$ for some prime $p$. Let $d$ be
the diagonal multiplier of $T$. We note that by
Lemma~\ref{centraliser} and Lemma~\ref{PSL-PSU-gcd} (and 
Hypotheses~\ref{HS-hyps}), we have
that $a/r\geq \gcd(|\Out T|,|T|-1)\geq 5$ if $T$ is not a Ree or a
Suzuki group and similarly $a\geq \gcd(|\Out T|, |T|-1)\geq 5$ if $T$
is  a Ree or a Suzuki group.

We use Table~$5$ of~\cite{Atlas} extensively in the rest of the proof.

\setcounter{case}{0}
\begin{case}
The simple group $T$ is untwisted and has rank $n$ at least $2$.
\end{case}

We first assume that $T$ is not $A_n(t)$. From
Table~$5$ in~\cite{Atlas}, we get $|\Out T| \leq 24a$. Recall from 
Hypotheses~\ref{HS-hyps} that $y=\gcd(|H:\Inn T|,|T|-1)\leq |\Out T|$.
Then, using
Lemmas~\ref{centraliser},~\ref{Lie-lower-bd} 
and~\ref{PSL-PSU-gcd} and the fact that $d \le 4$ for these groups
(see Table~$5$ in~\cite{Atlas}), for
$s\in S$ we have  
$$
(p^{a}-1)^n/4 \le |C_H(s)| \le 2(24a)(y+1)\le 2(24a)(a+1).\
$$
Since $y\geq 5$ and $3\nmid y$, this inequality is
satisfied only if $n=2$ and $p=2$. But if $n=2,p=2$, then $d=1$, $|\Out
T|\leq 6a$ and the inequality $$(2^a-1)^2\leq |C_H(s)|\leq
2(6a)(a+1),$$ is never satisfied.  

If $T=A_n(t)$, then  we have $d=\gcd(n+1,p^a-1)$ and $|\Out
T|=2da$. Now, using Lemmas~\ref{centraliser},~\ref{Lie-lower-bd} 
and~\ref{PSL-PSU-gcd}, for 
$s\in S$ we have  
$$\frac{(p^{a}-1)^n}{\gcd(n+1,p^a-1)} \le |C_H(s)| \le 4\gcd(n+1,p^a-1)a(a+1).$$
Since $a\geq 5$, this inequality is
never satisfied.  (This is straightforward to check by dividing into 
three cases: 
 $d \le p^a-1$ for $n \ge 4$; $d \le 4$, for $n=3$;
and finally $d \le 3$ for $n=2$, where if $p=2$ and $a$ is odd we have $d=1$.)

\begin{case}The group $T$ is a unitary group.
\end{case}

Set $T={^2A_n}(t)$. By~Table~$5$ in~\cite{Atlas}, we have 
$d=\gcd(n+1,p^{a/2}+1)$ and $|\Out T|=da$. Now, using
Lemmas~\ref{centraliser},~\ref{Lie-lower-bd}  
and~\ref{PSL-PSU-gcd}, for 
$s\in S$ we have  
$$\frac{(p^{a/2}-1)^n}{\gcd(n+1,p^{a/2}+1)} \le |C_H(s)| \le
2\gcd(n+1,p^{a/2}+1)a(a/2+1).$$ 
Since $a/2\geq 5$, this inequality is satisfied only if
$n=2,p=2,a=10$. But if $T={^2A_2(2^{10})}$, then $|T|\equiv 2\pmod 5$, and 
so $\gcd(|\Out T|,|T|-1)=1<5$, contradicting Lemma \ref{centraliser}. 
(The impossibility of satisfying this inequality with other values of $n$, $p$ and $a$ is
straightforward to check by breaking it down into three cases: 
$d \le p^a+1$ for $n \ge 4$; $d \le 4$ for $n=3$;
and finally $d \le 3$ for $n=2$.)

\begin{case}
The simple group $T$ is twisted with $r=2$.
\end{case}

Because of Case 2, we may assume that $T$ is not a unitary
group. If $T$ is not a Suzuki or a Ree group (i.e. $T$ is $^2D_n(t)$ or
$^2E_6(t)$), then $d\leq 4$ and $|\Out T|\leq 4a$. In particular,
by Lemmas~\ref{centraliser},~\ref{Lie-lower-bd} and~\ref{PSL-PSU-gcd}, we 
have $$\frac{(p^{a/2}-1)^n}{4}\leq |C_H(s)|\leq 2(4a)(a/2+1).$$
Since $a/2\geq 5$, this inequality is never satisfied.

If $T$ is a Suzuki group or a Ree group (i.e. $T$ is
${}^2B_2(2^a),{}^2G_2(3^a)$, or ${}^2F_4(2^a)$), then $d=1$, $a$ is odd and
$|\Out T|=a$. By Lemmas~\ref{centraliser},~\ref{Lie-lower-bd}
and~\ref{PSL-PSU-gcd}, we have
that $$(p^{a/2}-1)^n\leq |C_H(s)|\leq 2a(a+1).$$
Since $a\geq 5$, this inequality is satisfied only if $p=2,n=2$ and $a=5,7$
(i.e. $T={^2B_2(2^5)},{^2B_2(2^7)}$). If $T={^2B_2(2^5)}$, then $|\Out
T|=5$
divides the order of 
$T$, therefore
$\gcd(|\Out T|,|T|-1)=1<5$, contradicting Lemma \ref{centraliser}. If
$T={^2B_2(2^7)}$, 
then $|\Out T|=7$ and $|T|\equiv 6\mod 7$, therefore $\gcd(|\Out
T|,|T|-1)=1<5$, contradicting Lemma \ref{centraliser}.

\begin{case}
The simple group $T$ is the Steinberg triality group $^3D_4(t)$.
\end{case}

From Table~$5$ in~\cite{Atlas}, we have $d=1$ and $|\Out T|=a$. By
Lemmas~\ref{centraliser},~\ref{Lie-lower-bd} and~\ref{PSL-PSU-gcd}, we
have $$(p^{a/3}-1)^4\leq |C_H(s)|\leq 2a(a/3+1),$$ 
which is never satisfied for $a/3\geq 5$.

\begin{case}
The simple group $T$ is untwisted of rank 1, and is not a unitary group.
\end{case} 

This can only happen if $T$ is the projective special linear group
$A_1(t)$.  The order of the outer automorphism group of $T$ is $a
\gcd(2,p^a-1)$, which is $a$ if $p=2$, and $2a$ 
otherwise.  

Now using Lemmas~\ref{centraliser} and~\ref{Lie-lower-bd}, for any $s\in
S$, we have  
$$\frac{p^a-1}{d} \le |C_H(s)| \le \begin{cases}2a(a+1) \quad\text{ if $p$ is
    even,}\\ 4a(a+1)\quad\text{ if $p$ is odd,}\end{cases}$$ 
where $d=1$ if $p$ is even, and $d=2$ if $p$ is odd. Now $a\geq 5$, so
this inequality is only satisfied when $p^a=2^5$.  
 
 If $T=A_1(32)=\PSL(2,32)$, we have $|\Out T|=5$, so $|H:\Inn T| \in \{1, 5\}$.
 By Lemma~\ref{centraliser}, $|H:\Inn T|\geq 5$. Therefore $|H:\Inn
 T|=5$ and  $H=\PGamL(2,32)$. Again by
 Lemma~\ref{centraliser}, $|C_H(s)| \le 30$ if $S$ consists of a
 single $H$-orbit, and $|C_H(s)| \le 60$ if $S$ is the union of two
 $H$-orbits.  Now, Lemma~\ref{Lie-lower-bd} gives $|C_H(s)| \ge 31$,
 so $S$ must be the union of two $H$-orbits, and by
 Lemma~\ref{connection}, the orbits must have the form $B$ and
 $B^{-1}$.  There are 6 conjugacy classes of elements of $\PSL(2,32)$
 in $\PGamL(2,32)$ for which the centralisers of the elements have
 order at most 60.  First, there is a single conjugacy class of
 elements of order 11; since there is only one such class, it cannot
 be in $S$.  Next, there are three conjugacy classes of elements of
 order 31; each of these is the conjugacy class represented by a
 matrix of the
 form  
$$\begin{bmatrix} x & 0 \\ 0 & x^{-1}\end{bmatrix},$$ which is
 conjugate in $\PSL(2,32)$ to its own inverse, since they have the
 same eigenvalues.  Therefore, these three conjugacy classes
 are inverse-closed, so they cannot be in $S$.  Finally,
 there are two conjugacy classes of elements of order 33.  The
 normaliser of a cyclic group of order 33 is a dihedral group of order
 66, where the involution inverts the element in the cycle, so these
 conjugacy classes are inverse-closed, and again cannot be in $S$.  Thus, $T \neq \PSL(2,32)$.  
 \end{proof} 

We summarise the results of this section in the following corollary.
\begin{cor}\label{noHS}
Let $\Gamma$ be a connected ve-srg, with
$G \le \Aut(\Gamma)$ acting transitively on the vertices and the edges.  Then $G$ cannot be a
holomorphic simple group.
\end{cor}

\begin{proof}
This is an immediate consequence of Theorems \ref{alternating} and \ref{Lie-type}, 
together with the Classification of Finite Simple Groups.
\end{proof}

\end{document}